\input amstex
\documentstyle{amsppt}
\magnification 1200
\topmatter
\title
Point processes and \\ the infinite symmetric group.\\
Part VI: Summary of results
\endtitle
\author
Alexei Borodin and Grigori Olshanski
\endauthor
\thanks Supported by the Russian Foundation for Basic Research under
grant 98-01-00303 (G.~O)
and by the Russian Program for Support of Scientific Schools under
grant 96-15-96060 (A.~B. and G.~O.)
\endthanks
\abstract
We give a summary of the results from Parts I--V

(math/9804086, math/9804087, math/9804088, math/9810013, math/9810014).

Our work originated from harmonic analysis on the infinite symmetric
group. The problem of spectral decomposition for certain
representations of this group leads to a family of probability
measures on an infinite--dimensional simplex, which is a kind of dual object for the infinite symmetric group. 

To understand the nature of these measures we interpret them as stochastic point processes on the punctured real line and compute their correlation functions. 

The correlation functions are given by multidimensional integrals which can be expressed in terms of a multivariate hypergeometric series (the Lauricella function of type B).

It turns out that after a slight modification (`lifting') of the processes the correlation functions take a common in Random Matrix Theory (RMT) determinantal form with a certain kernel.

The kernel is expressed through the classical Whittaker functions. 
It depends on two parameters and admits a variety of degenerations.
They include the well--known in RMT sine and Bessel kernels as well as some other Bessel--type kernels which, to our best knowledge, are new.

The explicit knowledge of the correlation functions enables us to derive a number of conclusions about the initial
probability measures. 

We also study the structure of our kernel; this finally leads to a constructive description of the initial measures.  

We believe that this work provides a new promising connection between RMT and Representation Theory.

\endabstract
\endtopmatter
\define\pz{P_{zz'}}

\define\rz{\rho_n^{(zz')}}
\define\wrz{\widetilde\rho_n^{(zz')}}
\define\sz{\sigma_n^{(zz')}}
\define\al{\alpha}
\define\be{\beta}
\define\su{\sum\limits}
\define\inte{\int\limits}
\define\sc{\widetilde s_{\lambda}}
\define\wt{\widetilde}
\TagsOnRight
\NoRunningHeads
 In this paper we review our results from [Part I -- Part V].
\subhead 1. The Thoma simplex [Part I, \S1]\endsubhead
The starting point of our study was the work [KOV] on generalized regular
representations of the infinite symmetric group. These representations depend
on two parameters (further denoted by $z$ and $z'$). The decomposition of the
representations into irreducibles is governed by certain probability measures
$P_{zz'}$ living on the infinite--dimensional simplex
$$
\Omega=\{\al_1\geq\al_2\geq\ldots\geq 0;\,\be_1\geq\be_2\geq\ldots\geq
0\,|\,\su_{i=1}^{\infty}(\al_i+\be_i)\leq 1\}
$$
called the {\it Thoma simplex\/} \cite{VK}, \cite{KV}. Note that
$\Omega$ is compact in the topology of pointwise convergence.

Our aim is to understand these measures. Our results show that the measures
$\pz$ are close to stochastic point processes arising in scaling limits of
random matrix ensembles.

\subhead 2. Probability measures $\pz$ [Part I, \S\S1,2]\endsubhead
The measures $\pz$ are defined as follows. There exists a family
$\{\widetilde s_{\lambda}\}$ of continuous functions on $\Omega$, indexed by
the Young diagrams $\lambda$. Their linear span is dense in $C(\Omega)$. We
know explicitly the integrals $$
\varphi_{zz'}(\lambda)=\int_{\Omega}\widetilde
s_{\lambda}(\omega)\pz(d\omega)
\tag 2.1
$$
which determine $\pz$ uniquely. Essentially, this is the only information
about $\pz$ that we possess.

Now we shall describe the functions $\sc$ and write down the formula for
$\varphi_{zz'}$.

Let $\omega=(\al|\be)$ range over $\Omega$. Set
$$
\widetilde p_k(\omega)=\cases 1,& k=1\\
                              \su_{i=1}^{\infty}\al_i^k+(-1)^{k-1}
                              \su_{i=1}^{\infty}\be_i^k,& k\geq 2
                         \endcases
$$
and for any partition $\rho=(\rho_1,\ldots,\rho_l)$
$$
\widetilde p_\rho(\omega)=\widetilde p_{\rho_1}(\omega)\cdots \widetilde
p_{\rho_l}(\omega).
$$
These are continuous functions on $\Omega$.
The functions $\sc$ are related to $\wt p_\rho$'s in exactly the same way as
the Schur functions $s_\lambda$ are related to the power sums $p_\rho$, see
[M, \S1.7]. Specifically,
$$
\sc=\su_\rho \chi_\rho^\lambda z_\rho^{-1}\wt p_\rho
$$
where $\rho$ ranges over the partitions of
$|\lambda|=\lambda_1+\lambda_2+\ldots$; $\chi^\lambda$ is the irreducible
character of the symmetric group of degree $|\lambda|$, $\chi_\rho^\lambda$ is
its value on the conjugacy class indexed by $\rho$, and $z_\rho^{-1}\cdot
|\lambda|!$ is the cardinality of this conjugacy class.

 Conversely,
$$
\wt p_\rho=\su_\lambda\chi_\rho^\lambda\,\sc.
\tag 2.2
$$

The functions $\sc$ are called {\it extended Schur functions}, see [KV].

To define $\varphi_{zz'}(\lambda)$ we need Frobenius notation for a Young
diagram $\lambda$:
$$
\lambda=(p_1,\ldots,p_d|\, q_1,\ldots,q_d);
$$
here $d$ is the number of diagonal boxes of $\lambda$, and
$$
p_i=\lambda_i-i,\qquad q_i=\lambda_i'-i,
$$
where $\lambda'$ stands for the transposed diagram (see [M, \S1.1]).

Set $n=|\lambda|,\ t=zz'$. Then
$$
\aligned
\varphi_{zz'}(\lambda)=&\frac {t^d}{(t)_n}\prod\limits_{i=1}^d\frac{
(z+1)_{p_i}(z'+1)_{p_i}(-z+1)_{q_i}(-z'+1)_{q_i}}{(p_i!)^2(q_i!)^2}
\\ \times & \frac {\prod_{i<j}[(p_i-p_j)(q_i-q_j)]}{\prod_{i,j}(p_i+q_j+1)}
\endaligned
\tag 2.3
$$
where $(a)_m=a(a+1)\cdots(a+m-1)$, $(a)_0=1$.

We shall always assume that $z$ and $z'$ satisfy one of the following
conditions
$$
\aligned
&(*)\quad\ \, z'=\bar{z},\ z\in \Bbb C\setminus \Bbb Z;\\
&(**)\quad z\ \text{and}\ z'\ \text{are real and for a certain}\ m\in \Bbb Z,\
m<z,z'<m+1.
\endaligned
\tag 2.4
$$
Under this assumption there exists a probability measure $P_{zz'}$ such that
(2.1) holds.

Note that the measure $P_{zz'}$ can be obtained as the limit, as
$n\to\infty$, of certain statistics on partitions of $n$, see [KOO, proof of
Theorem B], [Part II, Introduction].

\subhead 3. Stochastic point processes [Part I, \S4]\endsubhead
We shall interpret the measures $\pz$ as stochastic point processes on
$I=[-1,1]\setminus\{0\}$. With each point $\omega\in \Omega$ we associate a
point configuration in $I$,
$$
\omega=(\al|\, \be)\mapsto (\al_1,\al_2,\ldots;-\beta_1,-\be_2,\ldots)
$$
where we omit possible zeros in $\al$ and $\be$. Note that any such
configuration has no accumulation points in $I$, the points condensate near
the origin which is not in $I$.

Thus, the measure $\pz$ becomes a probability measure on the space of point
configurations in $I$, i.e., a stochastic point process on $I$. We shall
denote this process by $\Cal P_{zz'}$.

The $n$th correlation function $\rho_n(x_1,\ldots,x_n)$ of a point process is
the density of the probability to find a point in each of the infinitesimal
intervals $[x_i,x_i+dx_i]$ (see [DVJ], [Part I, \S4] for details).

Our aim is to compute the correlation functions $\rz$ of our processes $\Cal
P_{zz'}$.

\subhead 4. Moment problems [Part I, \S3], [Part II, Section 1.1]\endsubhead
Our strategy to solve the problem is to find the probability measures $\sz$
on $[-1,1]^n$, $n=1,2,\ldots$, characterized by their moments
$$
\inte_{[-1,1]^n}x_1^{l_1}\cdots x_n^{l_n}\sz(dx)=\inte_{\Omega}\wt
p_{(l_1+1,\ldots,l_n+1)}(\omega)\pz(d\omega).
$$
Note that the RHS is explicitly known because of (2.1), (2.2):
$$
\inte_{[-1,1]^n}x_1^{l_1}\cdots x_n^{l_n}\sz(dx)=
\su_{\lambda,\,\,|\lambda|=\su_i(l_i+1)}
\chi^\lambda_{(l_1+1,\ldots,l_n+1)}\,\varphi_{zz'}(\lambda).
\tag 4.1
$$
It turns out that outside the diagonals $x_i=x_j$ in $I^n$
$$
\rz(x_1,\ldots,x_n)=\frac 1{|x_1\cdots x_n|}\, \{ \text{density of}
\ \sz(dx_1,\ldots,dx_n)\}
\tag 4.2
$$
It is worth noting that the measure $\sz$ for $n\geq 2$ always has singular
components living on the diagonals while (as is proved in [Part I,
Proposition 4.2], [Part II, Theorem 2.5.1]) $\rz$ has no such
components\footnote{This just means that nonzero coordinates in $\al$ and
$\be$ are pairwise distinct with probability 1.}.

Thus, we obtain the correlation functions from a more sophisticated object.
However, we can not formulate a moment problem for $\rz$ directly because of
the absolute value sign in (4.2).

\subhead{5. Integral represenations for the correlation functions
[Part I, \S\S 5,6], [Part II, Chapter 2]}\endsubhead The moment problem (4.1)
can be completely solved.  In particular, for $n=1$ we obtained a two
dimensional integral represenation of $\sigma_1^{(zz')}$ (see [Part I,
Theorem 5.8]). A careful examination shows that $\sigma_1^{(zz')}$ has no
atom at zero. This fact has an important corollary.  \proclaim{Theorem I
([Part I, Theorem 6.1])} With probability 1, $$
\su_{i=1}^{\infty}(\al_i+\be_i)=1.  $$ \endproclaim For $n>1$ the solution of
the moment problem (4.1) requires a lot of work ([Part II, Chapter 1]). The
difficulties come from the fact that the RHS of (4.1) involves the symmetric
group characters for which there is no simple expression. We use
Murnaghan--Nakayama rule ([M, \S1.7, Ex. 5]) to handle the characters. In
[Part II, Theorem 1.2.1'] we obtain a more economic form of this rule which
enable us to solve the moment problem.

The final expression for $\rz$ is a linear combination of multidimensional
integrals of various orders up to $3n$. The situation simplifies when all the
variables $x_1,\ldots, x_n$ are of the same sign, say, positive.

\proclaim{Theorem II ([Part II, Theorem 2.2.1])} Let $x_1,\ldots,x_n>0$. Then
$$
\gathered
\rz(x_1,\ldots,x_n)=t^n\Gamma(t)\\ \times \inte_{\Sb a_i,b_i>0\\
\su_ix_i(a_i+b_i+1)<1\endSb}\prod\limits_{i=1}^n\frac
{a_i^{-z}}{\Gamma(-z+1)} \frac {(a_i+1)^{z'}}{\Gamma(z'+1)} \frac
{b_i^{-z'}}{\Gamma(-z'+1)} \frac {(b_i+1)^{z}}{\Gamma(z+1)} \\ \times
\det\left(\frac 1{a_i+b_j+1}\right)\,
\frac{(1-\sum_ix_i(a_i+b_i+1))^{t-n-1}}{\Gamma(t-n)}\, dadb.
\endgathered
\tag 5.1
$$
\endproclaim
The RHS of (5.1) is well defined for $z$ and $z'$ such that
$$
-1<\Re z,\Re z'<1;\quad t=zz'>n.
$$
For other values of $z,z'$ we use analytic continuation.

We tacitely assume that $\sum_ix_i<1$. Actually, the correlation measure
$\rz(x)dx$ always lives on the set $\sum_i|x_i|\leq 1$; additional
considerations show that there are no singular components on the faces
$\sum_i|x_i|=1$, see the beginning of the proof of Theorem 3.3.1 in [Part
II].

The RHS of (5.1) can be expressed via the multivariate Lauricella
hypergeometric function of type B, see [Part II, Section 2.4]. In particular,
$\rho_1^{(zz')}$ can be expressed in terms of the Appell's two--dimensional
hypergeometric function $F_3$, see [Part II, Corollary 2.4.2]. Another
expression of $\rho_1^{(zz')}$ through the Lauricella function in three
variables is given in [Part I, Theorem 5.12].

\subhead 6. Lifting [Part II, Chapter 3]\endsubhead
A surprising fact is that the correlation functions are greatly simplified
after a `lifting' of our processes $\Cal P_{zz'}$ to a slightly bigger space.

According to Theorem I, the measure $\pz$ ic concentrated on the face
$$
\Omega_0=\{\omega=(\al\,|\, \be)\,|\, \su_i(\al_i+\be_i)=1\}\subset \Omega.
$$
Let
$$
\wt\Omega_0=\Omega_0\times \Bbb R_+.
$$
We pass from the measure $\pz$ on $\Omega_0$ to the measure
$$
\wt \pz =\pz\otimes \left\{\frac {s^{t-1}}{\Gamma(t)}\, e^{-s}ds\right\}
$$
on $\wt\Omega_0$. In other words, we tensor $\pz$ by the gamma--distribution
with parameter $t$.

We associate to a point $\wt\omega=(\omega,s)\in \wt \Omega_0$ a point
configuration in $\Bbb R^*=\Bbb R\setminus \{0\}$ as follows
$$
\aligned
((\al|\,\be),s)\mapsto&(\wt\al_1,\wt\al_2,\ldots;-\wt\be_1,-\wt\be_2,\ldots)\\
=&(s\al_1,s\al_2,\ldots;-s\be_1,-s\be_2,\ldots).
\endaligned
$$
Thus, we get a probability measure on the space of point configurations in
$\Bbb R^*$, i.e., a point process.
We denote it by $\wt{\Cal P_{zz'}}$.

The process $\wt{\Cal P_{zz'}}$ is obtained from $\Cal P_{zz'}$ by multiplying
the random configuration in $I$ by the scalar random factor $s$ with
gamma--distribution.

The lifting is invertible via the map
$$
(\wt \al|\, \wt\be)\mapsto (\al|\, \be)=\left(\frac{\wt\al_1}s,
\frac{\wt\al_2}s,\ldots;-\frac{\wt\be_1}s,-\frac{\wt\be_2}s,\ldots\right)
$$
where $s=\sum_i(\wt\al_i+\wt\be_i)$.

Note that $s$ is finite almost surely with respect to $\wt\pz$.

Let $\wrz(x_1,\ldots,x_n)$ be the $n$th correlation function of $\wt{\Cal
P_{zz'}}$.

The `lifted' correlation functions are related to the initial ones by the
following simple transformation
$$
\wrz(x_1,\ldots,x_n)=\inte_0^\infty\frac{s^{t-1}}{\Gamma(t)}\,
\rz(x_1s^{-1},\ldots,x_ns^{-1})\frac {ds}{s^n},
\tag 6.1
$$
see [Part II, Proposition 3.1.1]. This implies that the moments of the
correlation measures $\rz(x)dx$ and $\wrz(x)dx$ are related in a very
simple way
$$
\gather
\int x_1^{l_1}\cdots x_n^{l_n}\wrz(x_1,\ldots,x_n)dx_1\ldots dx_n\\
=(t)_{l_1+\ldots+l_n}
\int x_1^{l_1}\cdots x_n^{l_n}\rz(x_1,\ldots,x_n)dx_1\ldots dx_n
\endgather
$$
(the moments are finite if $l_1,\ldots,l_n\geq 1$).

Note also that the transformation (6.1) is invertible, and there exists a
complex inversion formula similar to that for the Laplace transform.

\subhead 7. Determinantal formula. Matrix Whittaker kernel [Part IV, Sections
1,2] \endsubhead Now we shall show that the correlation functions $\wrz$ have
determinantal form with a kernel expressed through the Whittaker function
$W_{\kappa,\mu}(x)$, $x>0$.

This function can be characterized as the only solution of the Whittaker
equation
$$
y''-\left(\frac 14 -\frac\kappa x+\frac {\mu^2-\frac 14}{x^2}\right)\, y=0
$$
such that $y\sim x^\kappa e^{-\frac x2}$ as $x\to+\infty$ (see [E1, Chapter
6]).  Here $\kappa$ and $\mu$ are complex parameters. Note that $$
W_{\kappa,\mu}=W_{\kappa,-\mu}.
$$
We shall employ the Whittaker function for real $\kappa$ and real or pure
imaginary $\mu$; then $W_{\kappa,\mu}$ is real.

Set
$$
a=\frac{z+z'}2,\quad \mu=\frac{z-z'}2,\quad \sigma=\sqrt{\sin(\pi z)\sin(\pi
z')}>0.
\tag 7.1
$$
It is often convenient to consider $(a,\mu)$ as new parameters of our
processes, instead of $(z,z')$. In terms of $(a,\mu)$ the basic assumptions
(2.4) take the following form
$$
\align
&\bullet\  a\  \text{is always real}\\
&\bullet\ \text{either } \mu\ \text{is pure imaginary}\\
&\bullet\ \text{or } \mu \text{ is real, } |\mu|<\frac 12, \text{ and there
exists } m\in \Bbb Z\text{ such that }
\\ &\quad\, m+|\mu|<a<m+1-|\mu| \\
&\bullet\  a \text{ is not an integer when } \mu=0.
\endalign
$$
Introduce the following functions on $\Bbb R_+$
$$
\matrix
&A_+(x)=\sqrt{x}\,W_{a+\frac 12,\mu}(x)\qquad
&A_-(x)=\sqrt{x}\,W_{-a+\frac 12,\mu}(x)\\
&B_+(x)=\sqrt{x}\,W_{a-\frac 12,\mu}(x)\qquad
&B_-(x)=\sqrt{x}\,W_{-a-\frac 12,\mu}(x)
\endmatrix
$$
\proclaim{Theorem III ([Part IV, Theorem 2.7])} The correlation functions of
the lifted process $\wt{\Cal P_{zz'}}$ have the form $$ \gathered
\wrz(u_1,\ldots,u_n)=\det\left[ K(u_i,u_j)\right]_{i,j=1}^n,
\\ n=1,2,\ldots;\quad u_1,\ldots,u_n\in \Bbb R^*,
\endgathered
\tag 7.2
$$
where the kernel $K(u,v)$ is conveniently written in the block form
$$
K(u,v)=\cases
K_{++}(u,v),\quad &u,v>0;\\
K_{+-}(u,-v),\quad &u>0\,,v<0;\\
K_{-+}(-u,v),\quad &u<0\,,v>0;\\
K_{--}(-u,-v),\quad &u,v<0;
\endcases
$$
with
$$
\aligned
K_{++}(x,y)=\frac
1{\Gamma(z)\Gamma(z')}\,&\frac{A_+(x)B_+(y)-B_+(x)A_+(y)}{x-y}\\
K_{+-}(x,y)=\frac
\sigma{\pi}\,&\frac{A_+(x)A_-(y)+tB_+(x)B_-(y)}{x+y}\\
K_{-+}(x,y)=-\frac
\sigma{\pi}\,&\frac{A_+(y)A_-(x)+tB_+(y)B_-(x)}{x+y}\\
K_{--}(x,y)=\frac
1{\Gamma(-z)\Gamma(-z')}\,&\frac{A_-(x)B_-(y)-B_-(x)A_-(y)}{x-y}
\endaligned
\tag 7.3
$$
(recall that $t=zz'$).
\endproclaim
This is one of our main results.

The matrix representation of the kernel
$$
K=\left[\matrix K_{++}&K_{+-}\\K_{-+}&K_{--}\endmatrix\right]
\tag 7.4
$$
correponds to the splitting $\Bbb R^*=\Bbb R_{+}\sqcup \Bbb R_{-}$ of the
phase space and subsequent identification of $\Bbb R_-$ with the second copy
of $\Bbb R_+$. We call (7.4) the {\it matrix Whittaker kernel}.

All the blocks of (7.4) are real valued kernels on $\Bbb R_+$. They possess
the following symmetry
$$
\gather
K_{++}(x,y)=K_{++}(y,x),\quad K_{--}(x,y)=K_{--}(y,x),\\
K_{+-}(x,y)=-K_{-+}(y,x).
\endgather
$$
Note the minus sign in the last relation. It means that the kernel (7.4) is
formally $J$--symmetric where $J$ is the operator
$\operatorname{id}\oplus\operatorname{(-id)}$ in
$L^2(\Bbb R_+,dx)\oplus L^2(\Bbb R_+,dx)$.

Another symmetry property: the change of parameters $(z,z')\to(-z,-z')$ is
equivalent to the transform of the kernel
$$
\left[\matrix K_{++}&K_{+-}\\K_{-+}&K_{--}\endmatrix\right]
\longrightarrow
\left[\matrix K_{--}&-K_{-+}\\K_{+-}&K_{++}\endmatrix\right]
$$

Determinantal form for the correlation functions (like (7.2)) appears in
different problems of random matrix theory and mathematical physics, see,
e.g., [Dy], [Me1], [Ma1], [Ma2], [TW1--3], [L], [KBI]. In most situations the
kernel $K$ is symmetric or Hermitian (see, however, [B]). Appearance of
$J$--symmetric kernels seems to be new.  

\subhead 8. The $L$--kernel [Part V, \S2]\endsubhead 
Consider the operator $K$ in the Hilbert space $$ L^2(\Bbb
R^*,du)\simeq L^2(\Bbb R_+,dx)\oplus L^2(\Bbb R_+,dx) \tag 8.1 $$ given by
the kernel (7.4).  

\proclaim{Theorem IV ([Part V, Theorem 2.4])} Assume $|a|<\frac
12$. Then 
$$ K=\frac{L}{1+L} \tag 8.2 $$ 
where $L$ is bounded and given by
the kernel 
$$ L(x,y)=\left[\matrix 0& {\frac \sigma\pi}\,\left(\frac
xy\right)^a\,\frac{\exp\left(-\frac{x+y}2\right)}{x+y}\\
-\frac\sigma\pi\,\left(\frac
yx\right)^a\,\frac{\exp\left(-\frac{x+y}2\right)}{x+y}& 0
\endmatrix
\right]
\tag 8.3
$$
Recall that $a$ and $\sigma$ were defined in (7.1).
\endproclaim

It is worth noting that (8.3) involves no special functions. Note also that
$\mu=\frac{z-z'}2$ enters only in the scalar factor $\sigma$.

The formulas (7.2), (7.3), (8.2), (8.3) give a precise description of the
process $\wt{\Cal P_{zz'}}$ and, thus, of the initial process $\Cal P_{zz'}$.
We consider these formulas as our main result.

\subhead 9. Spectral analysis [Part V, \S3]\endsubhead
In this section we shall diagonalize the operators $K$ and $L$.
Denote by $A$ the operator in $L^2(\Bbb R_+,dx)$ with the kernel
$$
A(x,y)=
\frac\sigma\pi\,\left(\frac
xy\right)^a\frac{\exp\left(-\frac{x+y}2\right)}{x+y}
\tag 9.1
$$
It is bounded provided that $|a|<\frac12$. Let $A'$ denote the operator in
$L^2(\Bbb R_+,dx)$ with the transposed kernel
$$
A'(x,y)=A(y,x).
$$
Then (8.1) implies that the blocks of the matrix Whittaker kernel are
expressed through $A$ and $A'$ as follows
$$
\matrix
K_{++}=AA'(1+AA')^{-1}\quad &
K_{--}=A'A(1+A'A)^{-1}\\
K_{+-}=(1+AA')^{-1}A\quad &
K_{-+}=-(1+AA')^{-1}A'
\endmatrix
\tag 9.2
$$
Consider the following Sturm--Liouville differential operator on $(\Bbb
R_+,dx)$ depending on $a$ as a parameter:
$$
\Cal D_a=\frac d{dx}\,x^2\frac d{dx} +\left(a-\frac x2\right)^2.
\tag 9.3
$$
We have
$$
\Cal D_aA=A\,\Cal D_{-a}
$$
which implies that $\Cal D_a$ (formally) commutes with $AA'$ and $K_{++}$
while $\Cal D_{-a}$ commutes with $A'A$ and $K_{--}$.

Consider the following functions on $\Bbb R_+$:
$$
f_{a,m}(x)=\frac 1x W_{a,im}(x),\quad m>0.
$$
We have
$$
\Cal D_a f_{a,m}=\left(a^2+\frac 14+m^2\right)f_{a,m}.
$$
According to [W], the functions $f_{a,m}$ with $a$ fixed form a continual
basis in $L^2(\Bbb R_+,dx)$ diagonalizing $\Cal D_a$. Moreover, the
Plancherel formula looks as follows. For good enough functions $f(x)$ and
$g(x)$
$$
(f,g)=\frac1{\pi^2}\inte_0^\infty(f,f_{a,m})(f_{a,m},g)\cdot \Gamma\left(\frac
12-a-im\right)\Gamma\left(\frac12-a-im\right)dm
$$
where $(\cdot,\cdot)$ stands for the inner product in $L^2(\Bbb R_+,dx)$.

We have
$$
Af_{a,m}=\frac
\sigma\pi\,\Gamma\left(\frac12-a+im\right)\Gamma\left(\frac12-a-
im\right)f_{a,m},
$$
$$
A'f_{a,m}=\frac \sigma\pi\,\Gamma\left(\frac12+a+im\right)\Gamma\left(\frac
12
+a-im\right)f_{-a,m}.
$$
Returning to the decomposition (8.1), we take $\{f_{a,m}\}_{m>0}$ as a basis
in the first summand $L^2(\Bbb R_+,dx)$ and $\{f_{-a,m}\}$ as a basis in the
second one.

Then we get a basis in the whole space $L^2(\Bbb R^*,du)$ diagonalizing both
$$
L=\left[\matrix 0&A\\-A'&0\endmatrix\right]\ \text{and}\ K=\left[\matrix
AA'(1+AA')^{-1}&
(1+AA')^{-1}A\\
-(1+A'A)^{-1}A'&
A'A(1+A'A)^{-1}
\endmatrix
\right].
$$
In particular, we get the diagonalization of $K_{++}$:
$$
K_{++}f_{a,m}=\frac{\cos(2\pi\mu)-\cos(2\pi a)}{\cos(2\pi\mu)+\cos(2\pi
im)}\,f_{a,m}.
$$

Note that for $a=0$ the integral transform inverse to
$f\mapsto\{(f,f_{a,m})\}_{m>0}$ is the Kontorovich--Lebedev transform, see
[E2].

The above results show that the operators $K$ and $L$ with $a$ fixed and
$\mu$ varying form a commutative family.

\subhead 10. Applications [Part III, Sections 2,5]\endsubhead
Now we shall give applications of the main results. These applications
concern our initial object --- the probability measures $\pz$ on $\Omega$.
\proclaim{Theorem V ([Part III, Theorem 5.1])} Consider $\al_1,\al_2,\dots$;
$\be_1,\be_2,\dots$ as random variables on the probability space
$(\Omega,\pz)$. Then with probability 1 there exist the limits
$$
\lim\limits_{k\to\infty}\al_k^{1/k}=
\lim\limits_{k\to\infty}\be_k^{1/k}=
\exp\left(-\frac{\pi\sin(2\pi\mu)}{2\mu\sigma^2}\right).
\tag 10.1
$$
\endproclaim
This result is a kind of the strong law of large numbers. Roughly speaking,
it means that both $\al_i$'s and $\be_i$'s decay as geometric progressions
with the same exponent. Similar situation occurs for the Poisson--Dirichlet
process, see [VS].

The same limit relation as (10.1) holds for $\wt\al_k$'s and $\wt\be_k$'s.

The proof of (10.1) is based on the examination of
 $\widetilde\rho_1^{(zz')}$
and
 $\widetilde\rho_2^{(zz')}$, and it incorporates Kingman's argument from [Ki,
Section 4.2].

Observe that the exponent of the decay does not change under the shifts
$(z,z')\to(z+N,z'+N)$, $N\in \Bbb Z$. Moreover, the whole process describing
the tails of the sequences $\{\al_i\}$ and $\{\be_i\}$ turns out to be
invariant with respect to these shifts. This periodicity is quite unexpected:
initial formulas (2.3) have no apparent periodicity.

Another application relies on the examination of
 $\widetilde\rho_1^{(zz')}$ alone.
\proclaim{Theorem VI ([Part III, Proposition 2.2])} Let
$$
|\al|=\su_{i=1}^\infty\al_i,\quad 
|\be|=\su_{i=1}^\infty\be_i,\quad
\psi(\cdot)=\frac{\Gamma'(\cdot)}{\Gamma(\cdot)}.
$$

Then
$$
\Bbb E\,|\alpha|=\sin(\pi z)\sin(\pi
z')\left[\frac{z-z'}{2\pi\sin(\pi(z-z'))}\,
\frac{z+z'-1}{zz'}+\frac{\psi(-z')-\psi(-z)}{\pi\sin(\pi(z-z'))}\right]
$$
and $\Bbb E\,|\be|$ is given by the same formula with $(z,z')$ replaced by
$(-z,-z')$.
\endproclaim

Here $\Bbb E$ is the symbol of expectation.

\subhead 11. Tail process [Part III, Sections 3,4], [Part V, \S4]\endsubhead
Here we study the asymptotic behavior of the process $\wt{\Cal P_{zz'}}$ near
the origin. The starting observation is that
$$
\widetilde\rho_1^{(zz')}(u)\sim \frac c{|u|},\quad u\to 0
\tag 11.1
$$
where
$$
c=\frac{2\mu\sigma^2}{\pi\sin(2\pi\mu)}.
\tag 11.2
$$
This asymptotic relation agrees with the rate of decay of $\al_i$'s and
$\be_i$'s, see (10.1).

As before, we identify the phase space $\Bbb R^*$ with the disjoint union of
two copies of $\Bbb R_+$ and then in each copy of $\Bbb R_+$ we make the
following change of variable
$$
x\mapsto \xi=-c\ln x.
\tag 11.3
$$
\proclaim{Theorem VII ([Part V, Theorem 4.1])}
Let $K'(\xi,\eta)$ denote the matrix Whittaker kernel on $\Bbb R$ obtained
from the matrix Whittaker kernel (7.4) by the change of variable (11.3). Then
there exists the limit
$$
\lim\limits_{M\to+\infty}K'(\xi+M,\eta+M)=\Cal K(\xi,\eta).
$$
Here $\Cal K(\xi,\eta)$ is a translation invariant matrix kernel with $\Cal
K(\xi,\xi)\equiv 1$,
$$
\Cal K(\xi,\eta)=\left[\matrix
 \Cal K_{++}(\xi,\eta)&
 \Cal K_{+-}(\xi,\eta)\\
 \Cal K_{-+}(\xi,\eta)&
 \Cal K_{--}(\xi,\eta)
\endmatrix\right]=
\left[\matrix
F(\xi-\eta)&G(\xi-\eta)\\-G(\eta-\xi)&F(\xi-\eta)\endmatrix\right]
$$
where
$$
F(\xi)=\frac BA\,\frac{\sinh(A\zeta)}{\sinh(B\zeta)}
\tag 11.4
$$
$$
G(\zeta)=\frac 1{2\mu\sigma}\,
\frac{(\sin(\pi\mu)\cos(\pi a))\cosh(A\zeta)+(\sin(\pi
a)\cos(\pi\mu))\sinh(A\zeta)}{\cosh(B\zeta)}
$$
and
the constants $A,B$ are as follows
$$
B=\frac{\pi\sin(2\pi\mu)}{4\mu\sigma^2}>0,
$$
$$
A=2\mu B=\frac{\pi\sin(2\pi\mu)}{2\sigma^2}.
$$
\endproclaim
We call the point process on $\Bbb R\sqcup \Bbb R$ with the correlation
functions given by the determinantal formula with the matrix kernel $\Cal
K(\xi,\eta)$ the {\it tail process} for $\wt{\Cal P_{zz'}}$. By the
construction, it describes the behavior of $\wt\al_i$'s and $\wt\be_i$'s with
large $i$ after the appropriate rescaling.

In particular, the kernel $\Cal K_{++}(\xi,\eta)$ on $\Bbb R$ describes the
tail of
 $\{\wt\al_i\}$
 alone, and $\Cal K_{--}(\xi,\eta)$ does the same for
$\{\wt\be_i\}$. Since $\Cal K_{++}=\Cal K_{--}$, the tail properties of
 $\{\wt\al_i\}$
and
 $\{\wt\be_i\}$ are identical.

 The same kernel $\Cal K_{++}(\xi,\eta)$ appears in the scaling limit of the
correlation functions for the unlifted process $\Cal P_{zz'}$ restricted to
$(0,1]$. However, this requires more sophisticated considerations, see
[Part II, Sections 4.2, 4.3].

Recall that $\mu$ is either real or pure imaginary. According to this, the
constant $A$ is also either real or pure imaginary. In the latter case the
hyperbolic sine in the numerator of (11.4) turns into the ordinary sine.

The next result is parallel to Theorem IV.

\proclaim{Theorem VIII ([Part V, Proposition 4.2])} Let $|a|<\frac 12$. Then
$$
\Cal K=\frac {\Cal L}{1+\Cal L}
$$
where $\Cal L$ is a bounded operator in $L^2(\Bbb R,d\xi)\oplus L^2(\Bbb
R,d\xi)$ with the kernel
$$
\Cal L(\xi,\eta)=\left[\matrix 0&\frac{\sigma}{2\pi}\,
\frac {\exp(-2aB(\xi-\eta))}{\cosh(B(\xi-\eta))}\\
-\frac\sigma{2\pi}\,\frac {\exp(-2aB(\eta-\xi))}{\cosh(B(\eta-\xi))}&0
\endmatrix\right]
\tag 11.5
$$
where $B$ is as above.
\endproclaim
Finally, note that all formulas of this section are invariant under the
shifts $(z,z')\to (z+N,z'+N)$, $N\in \Bbb N$ (or, equivalently,
$(a,\mu)\to(a+N,\mu)$), cf. Section 10.

\subhead 12. Formalism of fermion point processes 
[Part III, Section 1], [Part V, \S1] \endsubhead 
Here we shall discuss elementary general properties of the point
processes with determinantal correlation functions \cite{DVJ},
\cite{Ma1}, \cite{Ma2}.

Let $\goth X$ be a phase space with reference measure $dx$, $K(x,y)$ be a
kernel on $\goth X$, and $\Cal P$ be a point process on $\goth X$ with the
correlation functions
$$
\rho_n(x_1,\ldots,x_n)=\det\left[K(x_i,x_j)\right]_{i,j=1}^n.
$$

Let $\goth Y$ be a subset of $\goth X$ such that
$$
\int_{\goth Y}\rho_1(x)dx<\infty.
$$
Then the number of points in $\goth Y$ is finite with probability 1.
Moreover, the probability to find exactly $n$ points located in the
infinitesimal volumes $dx_1,\ldots,dx_n$ around points $x_1,\ldots,x_n$
equals
$$
\pi_n(x_1,\ldots,x_n)dx_1\ldots dx_n=\frac{\det\left[L_{\goth
Y}(x_i,x_j)\right]_{i,j=1}^n}{\det(1+L_{\goth Y})}\, dx_1\cdots dx_n.
\tag 12.1
$$
Here $L_{\goth Y}(x,y)$ is the kernel of the operator $L_{\goth Y}$ in
$L^2(\goth Y,dx)$ such that
$$
K_\goth Y=\frac{L_{\goth Y}}{1+L_\goth Y}
$$
where $K_\goth Y$ is the operator in $L^2(\goth Y,dx)$ with the kernel
$K(x,y)$ restricted to $\goth Y$.

Thus, in case of finite point configurations the operator $L=\frac K{1-K}$
has a clear probabilistic meaning --- it gives the distribution functions
$\pi_n$.

It is tempting to apply (12.1) in our case, when the operator $L$ has
especially simple form. Unfortunately, we can not do this for the whole space
because the point configurations are infinite. Of course, we can restrict
ourselves to an appropriate region $\goth Y$, but then the simple form of our
$L(x,y)$ will be lost.

Again, in the case of finite configurations, if the space is a disjoint union
of two pieces and $L$ is written in block form according to this splitting,
then vanishing of the diagonal blocks of $L$ (as in (8.3) and (11.5)) exactly
means that the random configurations has equally many points in each of the
pieces (see \cite{Part V, Proposition 1.7}). We do not know how to
interpret such vanishing in our situation. 

\subhead 13. Distribution of $\wt\al_1$ [Part III, Section 2]\endsubhead
Here we shall consider the lifted process $\wt{\Cal P_{zz'}}$ restricted to
$\Bbb R_+\subset \Bbb R^*$; it is governed by the `$++$' block of the matrix
Whittaker kernel, see (7.3). We call $K_{++}(x,y)$ the {\it Whittaker
kernel}.

Note that for any $\tau>0$ the intersection of the random configuration with
$\goth Y(\tau)=[\tau,+\infty)$ is finite. By (12.1)
$$
\operatorname{Prob}\{\wt\al_1<\tau\}=\frac 1{\det(1+L_{\goth
Y(\tau)})}=\det(1-K_{\goth Y(\tau)}).
\tag 13.1
$$
Here, following the notation of Section 12, the kernel of $K_{\goth Y(\tau)}$
is obtained by restricting $K_{++}(x,y)$ on $[\tau,+\infty)$. As was pointed
out by Tracy [T2], a modification of the argument from [TW3, Section V.B]
allows to express the Fredholm determinant (13.1) through the
Painlev\'e transcendent V. 

\subhead 14. Comparison with random matrices. Degenerations of the
Whittaker
kernel [Part III, Sections 1,6], [Part V, \S5]\endsubhead
There are a lot of similarities between the processes $\wt{\Cal P_{zz'}}$ and
point processes arising from random matrices. Random matrix theory leads to
various kernels: the sine kernel, the Bessel kernel, the Laguerre kernel,
etc. (see \cite{Me1}, \cite{TW1--3}) All of them have the form
$$
\frac{\varphi(x)\psi(y)-\psi(x)\varphi(y)}{x-y}
\tag 14.1
$$
and so does the Whittaker kernel $K_{++}(x,y)$, see (7.3). Here
$\varphi(\cdot)$ and $\psi(\cdot)$ are solutions of certain linear
second--order differential equations. The Fredholm determinants of many
kernels of this form are expressed through Painlev\'e transcendents,
see [TW3, Sections III and V]. Various kernels (14.1) restricted to
suitable intervals 
commute with Sturm--Liouville operators (see \cite{G}, \cite{Me1,
\S5.3}, \cite{TW1}, \cite{TW2}). The same is true for the
Whittaker kernel restricted to $[\tau,+\infty)$, see [Part III, Section 6].

The Whittaker kernel degenerates to the Laguerre kernel of order $N$ and
parameter $\al>-1$ if we formally substitute $z'=N$, $z=N+\alpha$, see [Part
III, Remark 2.4]. The Bessel kernel is the scaling limit of the Laguerre
kernel as $N\to+\infty$. Similar scaling procedure applied to the Whittaker
kernel leads to a two--parametric family of Bessel--type kernels, see [Part
V, Theorem 5.1].

Likewise, the stationary kernels from Theorem VII generalize the sine kernel:
if $\mu\to i\infty$, (11.4) tends to $\frac{\sin(\pi\zeta)}{\pi\zeta}$
irrespective to the behavior of $a$. The kernels $\Cal K_{++}(\xi,\eta)$ have
already appeared in [BCM, MCIN] in connection with so--called $q$--Hermite
ensemble.

As for the matrix Whittaker kernel, it has a similarity with matrix kernels
arising from two--matrix ensembles ([EM], [Me2], [MS]), see [Part IV, Section
3].
\bigskip
\noindent{\bf Acknowledgements.}
We would like to thank Craig~A.~Tracy for his letter [T2] and for drawing our attention to the works [MTW] and [T1].

\bigskip
\bigskip

{\smc A.~Borodin}: Department of Mathematics, The University of
Pennsylvania, Philadelphia, PA 19104-6395, U.S.A.  E-mail address:
{\tt borodine\@math.upenn.edu}
\newline\indent
{\smc G.~Olshanski}: Dobrushin Mathematics Laboratory, Institute for
Problems of Information Transmission, Bolshoy Karetny 19, 101447
Moscow GSP-4, RUSSIA.  E-mail address: {\tt olsh\@iitp.ru,
olsh\@glasnet.ru}

\end